\documentclass[12pt]{article}

\usepackage{amstext}
\usepackage{amsfonts}
\usepackage{latexsym}
\usepackage{epsfig}
\usepackage{verbatim}

\usepackage{amsmath, amsthm, amssymb, amscd, epsfig}
\usepackage[dvips]{color}

\title{The Characterization of Radial Subspaces of Besov- and Lizorkin-Triebel Spaces by Differences}
\author{Winfried Sickel, Leszek Skrzypczak, and Jan Vyb\'\i ral
\thanks{Research supported by the DFG Research Center {\sc Matheon}
``Mathematics for key technologies'' in Berlin.}}

\date{\today}

\newtheorem{T}{Theorem}
\newtheorem{Lem}{Lemma}
\newtheorem{Prop}{Proposition}
\newtheorem{Def}{Definition}
\newtheorem{Cor}{Corollary}
\newtheorem{Rem}{Remark}

\newcommand{\bt}{\begin{T}}
\newcommand{\bl}{\begin{Lem}}
\newcommand{\bp}{\begin{Prop}}
\newcommand{\bc}{\begin{Cor}}
\newcommand{\bd}{\begin{Def}}
\newcommand{\br}[2]{\begin{Rem}\label{#1}{\rm #2}}
\newcommand{\er}{ \end{Rem}}
\newcommand{\et}{\end{T}}
\newcommand{\el}{\end{Lem}}
\newcommand{\ep}{\end{Prop}}
\newcommand{\ec}{\end{Cor}}
\newcommand{\ed}{\end{Def}}

\newcommand{\be}{\begin{equation}}
\newcommand{\ee}{\end{equation}}
\newcommand{\beq}{\begin{eqnarray}}
\newcommand{\eeq}{\end{eqnarray}}
\newcommand{\beqq}{\begin{eqnarray*}}
\newcommand{\eeqq}{\end{eqnarray*}}

\def\lsim{\raisebox{-1ex}{$~\stackrel{\textstyle <}{\sim}~$}}

\def\gsim{\raisebox{-1ex}{$~\stackrel{\textstyle >}{\sim}~$}}
\def\tr{\mathop{\rm tr \,}\nolimits}
\def\ext{\mathop{\rm ext \,}\nolimits}

\newcommand{\bpr}{{\bf Proof.}\ }
\newcommand{\epr}{\hspace*{\fill}\rule{3mm}{3mm}\\}
\newcommand{\bspq}{B^s_{p,q}}

\newcommand{\fspq}{F^s_{p,q}}

\newcommand{\Rd}{{\mathbb R}^{d}}
\newcommand{\R}{{\mathbb R}}
\newcommand{\N}{{\mathbb N}}
\newcommand{\Z}{\mathbb Z}
\newcommand{\Zd}{\Z^{d}}
\newcommand{\C}{{\mathbb C}}
\newcommand{\cs}{{\mathcal{S}}}
\newcommand{\ca}{{\mathcal{A}}}

\newcommand{\cl}{{\mathcal{L}}}
\newcommand{\wt}{\widetilde}

\newcommand{\cf}{{\mathcal{F}}\,}
\newcommand{\cfi}{{\mathcal{F}}^{-1}\,}

\newlength{\fixboxwidth}
\begin{document}

\maketitle

\begin{abstract}
\noindent
We continue our earlier investigations (see \cite{SSV1, SS2})
of radial subspaces of Besov and Lizorkin-Triebel
spaces on $\R^d$. This time we study characterizations of these subspaces by differences.
\end{abstract}

%&&&&&&&&&&&&&&&&&&&&&&&&&&&&&&&&&&&&&&&&&&&&&&&&&&&&&
%&&&&&&&&&&&&&&&&&&&&&&&&&&&&&&&&&&&&&&&&&&&&&&&&&&&&

\section{Introduction}

%&&&&&&&&&&&&&&&&&&&&&&&&&&&&&&&&&&&&&&&&&&&&&&&&&&&&&
%&&&&&&&&&&&&&&&&&&&&&&&&&&&&&&&&&&&&&&&&&&&&&&&&&&&&&&

The study of radial subspaces of spaces of differentiable functions was initiated at the end
of the seventies by  papers of Strauss \cite{Strauss}, Coleman, Glazer and Martin \cite{CGM}, Berestycki and Lions
\cite{BeLi}  and Lions \cite{Lions}. It became clear,
that the symmetry of a function in a Sobolev-type space context implies several remarkable properties, 
which do not extend to arbitrary functions in those spaces.
For example, the {\em Radial Lemma} of Strauss states, that
every radial function $f \in H^1 (\Rd), d\ge 2$, is almost everywhere equal to a function
$\wt{f}$, continuous for $x\neq 0$, such that
\be\label{strauss}
|\wt{f} (x)| \le c \, |x|^{\frac{1-d}{2}}\, \| \,f \, |H^1(\Rd)\|\, ,\quad x\in\R^d\setminus\{0\},
\ee
where $c$ depends only on $d$. The main aim of our papers \cite{SS,SSV1,SS2} has been the study
of the decay and regularity
properties of radial functions (very much in the spirit of the Radial Lemma of Strauss)
in the more general framework of Besov and Lizorkin-Triebel spaces. 
Let us  refer also to the recent paper of  Cho and Ozawa \cite{CO}
for some extensions of the Radial Lemma.
The philosophy used in our earlier papers  consists in the fact that all the information
about a (locally integrable) radial function $f:\R^d\to\C$ is contained in its trace
\be\label{trace}
(\tr f)(t):= f(t,0,\dots,0),\quad t\in \R.
\ee
The corresponding extension operator is then given by
\be\label{exten}
(\ext g)(x)=g(|x|),\quad x\in\R^d \, ,
\ee
where $g$ is an even function on $\R$.
In \cite{SSV1} and \cite{SS2} we have investigated properties of the  operators $\tr$
and $\ext$ in the frame of both homogeneous and inhomogeneous Besov and Lizorkin-Triebel spaces. 
Here we have extensively used adapted atomic decompositions
characterizing radial subspaces.
Afterwards this has been applied to describe the interesting interplay between regularity and decay properties of radial functions
(or distributions) from these classes.
In case of Besov- and Lizorkin-Triebel spaces on $\Rd$ probably the most transparent
characterizations are given in terms of differences.
One could ask whether the proofs, given in \cite{SSV1,SS2}, would become more transparent by working with differences.
For that reason we shall derive a characterization by differences
for the traces of elements of $RF^s_{p,q} (\Rd)$ and $RB^s_{p,q} (\Rd)$.
Not only for technical reasons we shall restrict ourselves to spaces with values of $s$ strictly less than $1$, 
see the detailed comments in Subsection \ref{erg}.
It turns out that even in this simplified situation
(because of $s<1$ it will be enough to work with first order differences)
the characterizations we have found do not seem to be a convenient tool  to derive
those  properties of radial functions as stated in the Lemma of Strauss. 
However, under certain extra conditions on the parameters the spaces
$\tr(RF^s_{p,q} (\Rd))$ and $\tr (RB^s_{p,q} (\Rd))$ allow an interpretation as weighted 
Besov and Lizorkin-Triebel spaces, where the weight is given by $w(t):= |t|^{d-1}$, $t\in \R$, see \cite{SSV1}.
Hence, our characterizations of $\tr(RF^s_{p,q} (\Rd))$ and $\tr (RB^s_{p,q} (\Rd))$ in terms of differences 
may be understood as a first hint 
how complicated those characterizations of weighted Besov and Lizorkin-Triebel spaces may look like in case of 
weights with singularities. 
\\
The paper is organized as follows.
In Section \ref{differenz} we derive the announced characterization by differences of the spaces
$\tr(RF^s_{p,q} (\Rd))$.
The next section is used for establishing 
characterization by differences of $\tr(RB^s_{p,q} (\Rd))$.
Section \ref{sec3}
is devoted to a comparison of 
$\tr(RF^s_{p,q} (\Rd))$ and $\tr(RB^s_{p,q} (\Rd))$
with some weighted spaces on the real line.
Finally, in the Appendix at the end of the paper we collect some 
definitions and properties of the function spaces under consideration.

%&&&&&&&&&&&&&&&&&&&&&&&&&&&&&&&&&&&&&&&&&
%&&&&&&&&&&&&&&&&&&&&&&&&&&&&&&&&&&&&&&&&&

\subsection*{Notation}

%&&&&&&&&&&&&&&&&&&&&&&&&&&&&&&&&&&&&&&&&&&&
%&&&&&&&&&&&&&&&&&&&&&&&&&&&&&&&&&&&&&&&&&&&

As usual, $\N$ denotes the natural numbers, $\N_0:= \N \cup \{0\}$,  $\Z$ denotes the integers and $\R$ the real numbers.
For the  complex numbers we use the symbol $\C$, for the Euclidean $d$-space  we use $\Rd$
and $\Zd$ denotes the collection of all elements in $\Rd$ having integer components.
Many times we shall use the abbreviation
\be\label{eq-03}
\sigma_{p,q} (d) := d\, \max \Big(0, \, \frac 1p - 1, \, \frac 1q - 1 \Big)\,.
\ee
The symbol $\sigma_{d-1}$ stands for the usual $d-1$-dimensional Hausdorff measure in $\R^d$.\\
%and $\omega_{d-1}:=\sigma_{d-1}(\{z \in \Rd:|z|=1\})$ denotes the area of the unit sphere.
At very few places we shall need the Fourier transform $\cf$ as well as its inverse transformation
$\cfi$, always defined on the Schwartz space  $\cs' (\Rd)$ of tempered distributions.
\\
If $X$ and $Y$ are two quasi-Banach spaces, then the symbol  $X \hookrightarrow Y$
indicates that the embedding is continuous.
The set of all linear  and bounded operators
$T : X \to Y$, denoted by $\cl (X,Y)$, is  equipped with the standard quasi-norm.
As usual, the symbol  $c $ denotes positive constants
which depend only on the fixed parameters $s,p,q$ and probably on auxiliary functions,
 unless otherwise stated; its value  may vary from line to line.
Sometimes we  will use the symbols ``$ \lsim $''
and ``$ \gsim $'' instead of ``$ \le $'' and ``$ \ge $'', respectively. The meaning of $A \lsim B$ is given by: there exists a constant $c>0$ such that
 $A \le c \,B$. Similarly $\gsim$ is defined. The symbol
$A \asymp B$ will be used as an abbreviation of
$A \lsim B \lsim A$.\\
If $E$ denotes a space of functions on $\R^d$ then by $RE$ we mean the subset of radial functions in $E$
and we endow this subset  with the same quasi-norm as the original space.
Inhomogeneous Besov and Lizorkin-Triebel spaces are denoted by
$\bspq$ and $\fspq$, respectively.
Definitions as well as some references are given in the Appendix.

%&&&&&&&&&&&&&&&&&&&&&&&&&&&&&&&&&&&&&&&&&&&&&&&&&&&&&&&&&&&&&&&&&&&&&&&&&&&&&&&
%&&&&&&&&&&&&&&&&&&&&&&&&&&&&&&&&&&&&&&&&&&&&&&&&&&&&&&&&&&&&&&&&&&&&&&&&&&&&&&&

\section{The characterization of radial Lizorkin-Triebel spaces  by differences}
\label{differenz}

%&&&&&&&&&&&&&&&&&&&&&&&&&&&&&&&&&&&&&&&&&&&&&&&&&&&&&&&&&&&&&&&&&&&&&&&&&&&&&&&
%&&&&&&&&&&&&&&&&&&&&&&&&&&&&&&&&&&&&&&&&&&&&&&&&&&&&&&&&&&&&&&&&&&&&&&&&&&&&&&&

In case of first order Sobolev spaces $RW^1_p (\Rd)$ there is a simple characterization of the trace spaces, see \cite{SSV1}.
Let $d \ge 2$ and  $1 \le p< \infty$.
The  mapping $\tr $ is a linear isomorphism (with inverse $\ext$)
of $RW^1_p(\Rd)$ onto
the closure of  $RC^{\infty}_0 (\R)$ with respect to the norm
\be\label{eq-04}
\|\, g\, |L_p(\R,|t|^{d-1})\| + \|\, g'\, |L_p(\R,|t|^{d-1})\|\, .
\ee
As usual,
\be\label{eq-40}
\|\, g\, |L_p(\R,|t|^{d-1})\|:=
\Big(\int_{-\infty}^ \infty |g(t)|^p\,|t|^{d-1}\, dt\Big)^{1/p}\, .
\ee
The characterization of $\tr (RW^1_p (\Rd))$ makes clear what type of result can be expected
in the general situation of Lizorkin-Triebel spaces:
we shall look for a characterization as a weighted Lizorkin-Triebel space on $\R$,
where the weight is given by $w(t):= |t|^{d-1}$, $t \in \R$.
However, the outcome will be a bit more technical.

%&&&&&&&&&&&&&&&&&&&&&&&&&&&&&&&&&&&&
%&&&&&&&&&&&&&&&&&&&&&&&&&&&&&&&&&&&&&

\subsection{Some preliminaries}

%&&&&&&&&&&&&&&&&&&&&&&&&&&&&&&&&&&&&&&
%&&&&&&&&&&&&&&&&&&&&&&&&&&&&&&&&&&&&&&

Now we look for a counterpart of (\ref{eq-04}) in case of the radial Lizorkin-Triebel spaces $RF^s_{p,q} (\Rd)$.
Recall in this context that $W^ 1_p (\Rd) = F^1_{p,2} (\Rd)$, $1 <p< \infty$, in the sense of equivalent norms.
Our point of departure is the following characterization of
$F^s_{p,q}(\Rd)$, for which we refer to \cite[Thm.~3.5.3]{Tr92}.
We shall use the abbreviations
\[
M_{t,u} f (x):= \Big(t^{-d}\, \int_{|h|<t} |f(x+h)-f(x)|^u \, dh \Big)^{1/u} \, , \qquad t>0, \quad 0 < u < \infty\, , 
\]
and 
\[
M_{t,\infty} f (x):=  \sup_{|h|<t}\,  |f(x+h)-f(x)| \, , \qquad t>0 \, . 
\]
Suppose $0 <p <\infty$, $0 < q \le \infty$, $1 \le v \le \infty$ and
\be\label{eq-01}
s > d\, \max \Big(0,  \frac 1p - \frac 1v\, ,~ \frac 1q - \frac 1v\Big)\, .
\ee
Let $0 < u \le v$,  $s < 1$ and $T >0$. Then $F^s_{p,q}(\Rd)$ is the collection of all
$f \in L_{\max(p,v)} (\Rd)$ s.t.
\be\label{eq-02}
\| \, f |F^s_{p,q} (\Rd)\|^* :=
\| \, f |L_p (\Rd)\| + \Big\| \Big(\int^T_0 t^{-sq} \, (M_{t,u} f(\, \cdot \,))^ q \, \frac{dt}{t}  \Big)^{1/q} \, \Big|L_p (\Rd)\Big\|\, .
\ee
Moreover, $\| \, \cdot \,  |F^s_{p,q} (\Rd)\|^*$ is  equivalent to
$\| \, \cdot \,  |F^s_{p,q} (\Rd)\|$
on $L_{\max(p,v)} (\Rd)$.

\begin{Rem}
 \rm
(i) By taking $v=1$ it becomes obvious that we have a characterization
with first order differences as long as $\sigma_{p,q} (d) < s < 1$, see (\ref{eq-03}).
\\
(ii) There are many references dealing with the characterization
of Lizorkin-Triebel spaces by differences (partly in a different form than here),
we refer in particular to Strichartz \cite{Str} ($1 <p< \infty, q=2$), Seeger
\cite{Se} and Triebel \cite[2.5]{Tr83}.
\end{Rem}

We need a modification of the above characterization. Let $0 < A < B < \infty$ be fixed real numbers.
Let $\Omega_t (x)$, $t>0$, $x\in \Rd$, be a family of open sets in $\Rd$ s.t.
\be\label{eq-05}
\{h \in \Rd: \: |h|< A \, t\} \quad \subset \: \Omega_t (x) \: \subset  \quad \{h \in \Rd: \: |h|<B\, t\} \, .
\ee
We define
\[
M^\Omega_{t,u} f (x):= \Big(t^{-d}\, \int_{\Omega_t (x)} |f(x+h)-f(x)|^u \, dh \Big)^{1/u} \, , \qquad t>0, \quad 0 < u \le \infty\, .
\]
Then we have the following.

\begin{Lem}\label{prep}
Suppose $0 <p <\infty$, $0 < q \le \infty$, $1 \le v \le \infty$,
$0 < u \le v$,
\be\label{eq-20}
d\, \max \Big(0\, , \frac 1p - \frac 1v \, ,  \frac 1q - \frac 1v  \Big)< s < 1
\ee
and  $T >0$.
Let $\Omega_t (x)$, $t>0$, $x\in \Rd$, be a family of open sets in $\Rd$ s.t.
(\ref{eq-05}) is satisfied for some $0 < A < B < \infty$.
Then $F^s_{p,q}(\Rd)$ is the collection of all
$f \in L_{\max(p,v)} (\Rd)$ s.t.
\be\label{eq-14}
\| \, f |F^s_{p,q} (\Rd)\|^* :=
\| \, f |L_p (\Rd)\| + \Big\| \Big(\int^T_0 t^{-sq} \, (M^\Omega_{t,u} f(\, \cdot \,))^ q \, \frac{dt}{t}  \Big)^{1/q} \, \Big|L_p (\Rd)\Big\|\, .
\ee
Moreover, $\| \, \cdot \,  |F^s_{p,q} (\Rd)\|^*$ is  equivalent to
$\| \, \cdot \,  |F^s_{p,q} (\Rd)\|$
on $L_{\max(p,v)} (\Rd)$.
\end{Lem}

\noindent
\bpr
As a first step one can prove the lemma in case
$\Omega_t (x) = \{h \in \Rd: \quad |h| \le C \, t\}$ for a fixed $C>0$
by following the original proof in \cite[3.5.3]{Tr92}.
Afterwards the claim follows by an obvious monotonicity argument.
\epr

%&&&&&&&&&&&&&&&&&&&&&&&&&&&&&&&&&&&&&&&&&&&&&&&&&&&&&&&&&&&&&&&&&&&&&&&&&&&&&&&
%&&&&&&&&&&&&&&&&&&&&&&&&&&&&&&&&&&&&&&&&&&&&&&&&&&&&&&&&&&&&&&&&&&&&&&&&&&&&&&&

\subsection{Differences and radial Lizorkin-Triebel spaces}

%&&&&&&&&&&&&&&&&&&&&&&&&&&&&&&&&&&&&&&&&&&&&&&&&&&&&&&&&&&&&&&&&&&&&&&&&&&&&&&&
%&&&&&&&&&&&&&&&&&&&&&&&&&&&&&&&&&&&&&&&&&&&&&&&&&&&&&&&&&&&&&&&&&&&&&&&&&&&&&&&

We start with an elementary observation.
Let $f$ be a radial function in $L_p (\Rd)$. Then the mapping $\tr$, see (\ref{trace}),
is an isomorphism onto the space $RL_p(\R, |t|^{d-1})$ with inverse $\ext$. Hence, 
$\tr f$ is well-defined for any $f \in RF^s_{p,q}(\R^d)$ s.t. $s > \sigma_{p,p} (d)$.
This is a consequence of the Sobolev type  embedding
\[
F^s_{p,q}(\R^d) \hookrightarrow L_u (\Rd)\, , \qquad p \le u \le \left\{ \begin{array}{lll}
\frac{d}{\frac dp -s} & \qquad & \mbox{if} \quad s < d/p;
\\
\infty && \mbox{if} \quad s> d/p\, ,
\end{array}\right.
\]
see, e.g., \cite[2.7.1]{Tr83}.
The main aim of this section is to establish characterizations by differences of radial subspaces of Lizorkin-Triebel spaces
by using $\tr f$ instead of $f$ itself.
To begin with we concentrate on the use of the means $M^\Omega_{t,\infty} $.

\begin{T}\label{thm:1}
Let $d  \ge 2$,  $0< p < \infty$, $0 < q \le \infty$ and 
\be\label{eq-20bb}
d\, \max \Big (\frac 1p \, , ~\frac 1q \Big) < s < 1\, .
\ee
Then the radial function $f \in L_p (\Rd)$ belongs to $F^s_{p,q}(\R^d)$ if, and only if,
$g :=\tr f$ satisfies
\beqq
%&& \hspace*{-1.0cm}
\| \, g \, \|^\# & := & \|\, g\, |L_p(\R,|t|^{d-1})\|
\\
& + &
\Big(\int_{-\infty}^\infty |r|^{d-1}
\Big[\int^1_0 t^{-sq} \, \sup_{-t \le w\le t} |g(r+ w) - g(r)|^q  \, \frac{dt}{t}  \Big]^{p/q} \, dr
\Big)^{1/p} < \infty \, .
\eeqq
Moreover, $\| \, g \, \|^\#$ is equivalent to $\|\, f \, |RF^s_{p,q}(\R^d)\|$.
\end{T}

\noindent
\bpr
{\em Step 1.} Preparations. 
By $\langle x,y\rangle$ we denote the scalar product of $x,y \in \Rd$.
We are going to use Lemma \ref{prep} with a particular family $\Omega_t (x)$.
In case $t=1$ we put
\be\label{eq-07}
\Omega_1 (x) := \{h \in \Rd: \quad |x+h|\le 2 \} \qquad \mbox{if} \qquad  |x|\le 1
\ee
and
\beq\label{eq-08}
\Omega_1(x) & := &  \Big\{ h \in \Rd: \quad \exists \tau >0 \quad \exists y \in \Rd  \quad \text{s.t.} \quad x+h = \tau \, y\, ,
\nonumber
\\
 |y|& = &|x|, \quad \langle x,y \rangle > |x|^2- \frac12\ \text{and}\ |x|- \frac12< \tau |x|<|x|+ \frac 12\Big\}
\eeq
for $|x|>1$.

\begin{center}
\resizebox{5cm}{!}{\begin{picture}(0,0)%
\includegraphics{picture2.pstex}%
\end{picture}%
\setlength{\unitlength}{4144sp}%
\begingroup\makeatletter\ifx\SetFigFont\undefined%
\gdef\SetFigFont#1#2#3#4#5{%
  \reset@font\fontsize{#1}{#2pt}%
  \fontfamily{#3}\fontseries{#4}\fontshape{#5}%
  \selectfont}%
\fi\endgroup%
\begin{picture}(8124,8124)(439,-7723)
\put(6121,-4111){\makebox(0,0)[lb]{\smash{{\SetFigFont{34}{40.8}{\familydefault}{\mddefault}{\updefault}{\color[rgb]{0,0,0}$e_1$}%
}}}}
\put(5041,-3391){\makebox(0,0)[lb]{\smash{{\SetFigFont{34}{40.8}{\familydefault}{\mddefault}{\updefault}{\color[rgb]{0,0,0}$x$}%
}}}}
\end{picture}%
}\qquad\qquad\qquad
\resizebox{5cm}{!}{\begin{picture}(0,0)%
\includegraphics{picture4.pstex}%
\end{picture}%
\setlength{\unitlength}{4144sp}%
\begingroup\makeatletter\ifx\SetFigFont\undefined%
\gdef\SetFigFont#1#2#3#4#5{%
  \reset@font\fontsize{#1}{#2pt}%
  \fontfamily{#3}\fontseries{#4}\fontshape{#5}%
  \selectfont}%
\fi\endgroup%
\begin{picture}(4299,4524)(214,-5923)
\put(1711,-4021){\makebox(0,0)[lb]{\smash{{\SetFigFont{17}{20.4}{\familydefault}{\mddefault}{\updefault}{\color[rgb]{0,0,0}$e_1$}%
}}}}
\put(2746,-3976){\makebox(0,0)[lb]{\smash{{\SetFigFont{17}{20.4}{\familydefault}{\mddefault}{\updefault}{\color[rgb]{0,0,0}$x=2e_1$}%
}}}}
\end{picture}%
}\\
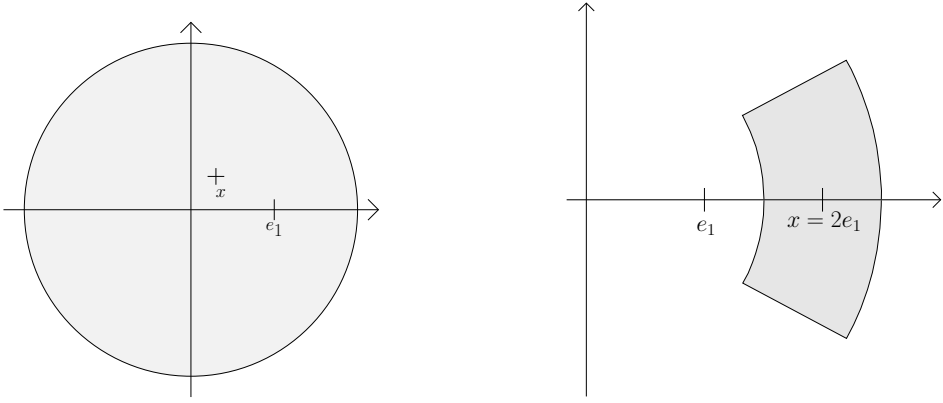
The set $x+\Omega_1(x)$ for $|x|\le 1$ (left) and $|x|>1$ (right).
\end{center}

For the general definition we make use of scaling, i.e., we define
\be\label{eq-09}
h \in \Omega_1 (x) \qquad \Longleftrightarrow \qquad t \, h \in \Omega_t (x)\, , \qquad t >0\, .
\ee
A scaling argument yields that the family $(\Omega_t (x))_{t,x}$ satisfies (\ref{eq-05}) if, and only if,
the subfamily  $(\Omega_1 (x))_{1,x}$ satisfies (\ref{eq-05}) (with $t=1$).
Now we investigate the subfamily with $t=1$.
The condition (\ref{eq-05}) is obviously satisfied in case  $|x|\le 1$. To prove it also for
$|x|>1$, we argue as follows. Let $h\in \Omega_1 (x)$. Clearly, $1/2 < \tau < 3/2$ and
\beqq
|h| & =  & |\tau \, y - x | =|\tau \, x - x + \tau\,  y- \tau \, x| \le |x|\cdot|\tau-1|+ \tau |x-y|
\\
& \le & \frac 12+\frac32\cdot|x-y|\le 2,
\eeqq
where the last step in the estimate follows from
\[
|x-y|^2=2|x|^2-2\, \langle x,y\rangle \le 2|x|^2-2|x|^2+2\cdot 1/2=1.
\]
This proves  $\Omega_1 (x) \subset \{h\in \Rd:\quad |h|\le 2\}$.
Next we shall prove the relation $\{h\in \Rd:\quad |h|\le 1/4\} \subset \Omega_1 (x)$.
Let $|h|\le 1/4$.
We define
\[
\tau := \frac{|x+h|}{|x|} \qquad \mbox{and} \qquad  y:= \frac{|x|}{|x+h|} \cdot \, (x+h)\,  .
\]
Then $|y|=|x|$ and
$$
|x|-1/2< \tau \, |x|<|x|+1/2
$$
follows easily. Finally, we have to show, that the conditions $|x|\ge 1$ and $|h|\le 1/4$ imply
\begin{equation}\label{eq:equiv2}
\left\langle x,\frac{x+h}{|x+h|} \right\rangle \cdot |x|\ge |x|^2- \frac 12 \, ,
\end{equation}
see (\ref{eq-08}).
First we claim
\be\label{eq-12}
\langle x,x+h \rangle \ge |x+h|\cdot \sqrt{|x|^2-|h|^2}\,.
\ee
But
\beqq
\Big(\langle x,x+h \rangle \Big)^2  & = &  |x|^4 + 2\, \langle x,h \rangle \, |x|^2
+ \langle x,h \rangle ^ 2
\\
& \ge &
|x|^ 4 + 2\, \langle x,h \rangle\, |x|^2 + |h|^2 \, |x|^2 - |h|^2\, |x|^2 -|h|^4
- 2\, \langle x,h \rangle\, |h|^2
\\
& = & \langle x+h, \, x+h \rangle\, (|x|^2 - |h|^2)\, .
\eeqq
Since $\langle x,x+h \rangle >0$ the claim (\ref{eq-12}) follows.
Hence
\[
\langle x, x+h \rangle \cdot \frac{|x|}{|x+h|} \ge
|x| \cdot \sqrt{|x|^2-|h|^2} \ge |x|^2- \frac 12 \, ,
\]
since $|x|\ge 1$.
This proves that our family $(\Omega_1 (x))_x$ and therefore
$(\Omega_t (x))_{t,x}$ satisfies (\ref{eq-05}) with $A=1/4$ and $B=3$.
\\
{\em Step 2.}
Let $\ca:\R^d\to\R^d$ be a rotation around the origin.
Hence, $\ca$ is a $(d\times d)$ matrix. The mapping $x \mapsto \ca (x)$ is linear,
length and angles are invariant under this mapping.
Simple calculations show, that $\Omega_t (\ca (x)) = \ca (\Omega_t(x))$.
\\
Now, let $x,z$ be two points in $\Rd$ s.t. $|x|=|z|$ and let $\ca$ be any rotation
around the origin satisfying $\ca (x) = z$.
For the moment we shall deal with $M^\Omega_{t,u} $ for all $u$, $0 < u \le \infty$.
Let  $f \in RF^s_{p,q} (\Rd)$. With  $g:= \tr f$ we conclude
\begin{eqnarray*}
M^\Omega_{t,u} f (z) & = & \Big( t^{-d}\,
\int_{\Omega_t (z)} | g(|z+h|) - g(|z|)|^u \, dh
\Big)^{1/u}
\\
& = &
\Big( t^{-d}\,
\int_{\Omega_t (\ca (x))} | g(|\ca (x)+h|) - g(|\ca (x)|)|^u \, dh
\Big)^{1/u} \, .
\end{eqnarray*}
A transformation of coordinates $h':= \ca^{-1} (h)$
and $|\det \ca^{-1}| = 1 $ yield
\begin{eqnarray*}
M^\Omega_{t,u} f (z)
&= & \Big( t^{-d} \int_{\Omega_t(x)} |g(|\ca (x+h')|)-g(|\ca (x)|)|^u dh'\Big)^{1/u}
\\
&= & \Big( t^{-d} \int_{\Omega_t(x)} |g(|x+h'|)-g(|x|)|^u dh'\Big)^{1/u}
\\
&= & M^\Omega_{t,u} f (x) \, .
\end{eqnarray*}
Hence, $M^\Omega_{t,u} f $  is a radial function.
\\
{\em Step 3.} Let $e_1 \in \R^d$ be the unit vector in the first direction.
We employ Lemma \ref{prep}, Step 1  and Step 2. For the radial function $f \in RF^s_{p,q} (\Rd)$ and $g:= \tr f$ this yields
\beq
\label{eq-16}
\| \, f \, |F^s_{p,q} (\Rd)\| & \asymp &
\|\, g\, |L_p(\R,|t|^{d-1})\|
\nonumber
\\
&& +
\Big(\int_{0}^\infty r^{d-1}
\Big[\int^T_0 t^{-sq} \, (M^\Omega_{t,u} f(r e_1))^ q \, \frac{dt}{t}  \Big]^{p/q} \, dr
\Big)^{1/p} \, .
\eeq
From now we concentrate on $u=\infty$.
This causes the restrictions for $s$  in (\ref{eq-20bb}).
Then
\[
M^\Omega_{t,\infty} f (re_1) = \sup_{h \in \Omega_t(re_1)}
\, |f(re_1+h)-f(re_1)| \quad \left\{
\begin{array}{lll}
 & \le &  \sup\limits_{|w|\le 3t} \, |g(r+w)-g(r)|
\\
&&
\\
& \ge &  \sup\limits_{|w| \le t/4} \, |g(r+w)-g(r)|
\end{array}
\right.
\]
by using (\ref{eq-05}) and Step 1.
Consequently,
\beqq
&& \hspace*{-1.0cm}\|\, g\, |L_p(\R,|t|^{d-1})\|   +
\Big(\int_0^\infty r^{d-1}
\Big[\int^T_0 t^{-sq} \, \Big(\sup\limits_{|w| \le t/4} \, |g(r+w)-g(r)|\Big)^q \, \frac{dt}{t}  \Big]^{p/q} \, dr
\Big)^{1/p}
\\
&& \\
& \lsim &
\| \, f \, |F^s_{p,q} (\Rd)\|
 \lsim \|\, g\, |L_p(\R,|t|^{d-1})\|
\\
&& \\
& + &
\Big(\int_0^\infty r^{d-1}
\Big[\int^{T'}_0 t^{-sq} \, \Big(\sup\limits_{|w| \le 3t} \, |g(r+w)-g(r)|\Big)^q \, \frac{dt}{t}  \Big]^{p/q} \, dr
\Big)^{1/p}\, , 
\eeqq
where we can choose $T,T'>0$ as we want.
Using Lemma \ref{prep} with $T=1/3$ (estimate from below) and afterwards with $T=4$ (estimate from above) we conclude
\beqq
\| \, f \, |F^s_{p,q} (\Rd)\| & \asymp &
\|\, g\, |L_p(\R,|t|^{d-1})\|
\\
&& +
\Big(\int_0^\infty r^{d-1}
\Big[\int^1_0 t^{-sq} \, \Big(\sup\limits_{|w| \le t} \, |g(r+w)-g(r)|
\Big)^q \, \frac{dt}{t}  \Big]^{p/q} \, dr
\Big)^{1/p} \, .
\eeqq
Since $g$ is even we may replace $\int_0^\infty r^{d-1} \ldots \, dr$ by 
$\int_{-\infty}^\infty |r|^{d-1} \ldots \, dr$. The proof is complete.
\epr

\begin{Rem}\label{f-fall}
 \rm
The scale of Lizorkin-Triebel spaces generalizes three well-known scales of function spaces, namely
Sobolev spaces, Bessel potential spaces and Slobodeckij spaces.
In fact, we have
\begin{itemize}
 \item $W^m_p (\Rd) = F^m_{p,2} (\Rd)$, $1 <p< \infty$, $m \in \N_0$;
\item $H^s_p (\Rd) = F^s_{p,2} (\Rd)$, $1 <p< \infty$, $s \in \R$;
\item $W^s_p (\Rd) = F^s_{p,p} (\Rd)$, $1 \le p < \infty$, $s >0$, $s\not \in \N$,
\end{itemize}
see, e.g.,  \cite[2.2]{Tr83}.
Hence, Thm. \ref{thm:1} also gives characterizations of
$RH^s_p (\R)$, $1 <p< \infty$, $1/2 < s < 1$, and of  $RW^s_p (\Rd)$, $1 < p < \infty$, $d/p < s <1 $.
 \end{Rem}

Of course, the restrictions in $s$, given in (\ref{eq-20bb}), are rather inconvenient.
They are mainly caused by the use of the means $M^\Omega_{t,\infty} $.
For this reason we  turn now to the use of $M^\Omega_{t,u} $, $0 < u< \infty$.
Employing  (\ref{eq-05}) we find
\beqq%\label{mean1}
M^\Omega_{t,u} f(re_1) & \le  & \Big(t^{-d} \int_{|h|\le 3t} |f(re_1 + h) - f(re_1)|^u\, dh \Big)^{1/u}
\nonumber
\\
& = &
\Big(t^{-d} \int_{|w-re_1|\le 3t} |g(|w|) - g(r)|^u\, dw \Big)^{1/u}
\nonumber
\\
& = &
\Big(t^{-d} \int_{\max(0,r-3t)}^{r+3t}     |g(\lambda) - g(r)|^u\,
\sigma_{d-1} (Q_{\lambda,3t} (r) )\,      d\lambda \Big)^{1/u}
\eeqq
where
\be\label{qlambda}
Q_{\lambda,t} (r) := \{w\in \Rd:\: |w|= \lambda\, , \quad  |w-re_1|\le t\}\, .
\ee
Similarly, also by using (\ref{eq-05}),  we obtain some sort of reverse inequality
\[%\label{mean2}
M^\Omega_{t,u} f(re_1) \ge  \Big(t^{-d} \int_{\max(0,r-t/4)}^{r+t/4}     |g(\lambda) - g(r)|^u\,
\sigma_{d-1} (Q_{\lambda,t/4} (r) )\,      d\lambda \Big)^{1/u} \, .
\]
By obvious manipulations with $T$ this proves the following lemma.

\begin{Lem}\label{diff1}
Let $d \ge 2$, $0 <p< \infty$, $0 < q \le \infty$, $1 \le v < \infty$, $0 < u \le v$  and
\be\label{eq-20bc}
d\, \max \Big (0\, , \frac 1p - \frac 1v \, , ~\frac 1q - \frac 1v \Big) < s < 1\, .
\ee
Let $T >0$.
Then the radial function $f \in L_p (\Rd)$ belongs to $F^s_{p,q}(\R^d)$ if, and only if,
$g :=\tr f$ satisfies
\beq\label{eq-161}
&& \hspace*{-0.7cm}
\| \, g \, \|^\triangle := \|\, g\, |L_p(\R,|t|^{d-1})\|
\\
\nonumber
& & \hspace{-0.7cm}
\Big(\int_{0}^\infty r^{d-1} \Big[\int^T_0 t^{-sq}  \Big(t^{-d} \int_{\max(0,r-t)}^{r+t}     |g(\lambda) - g(r)|^u
\sigma_{d-1} (Q_{\lambda,t} (r) )  d\lambda \Big)^{\frac qu} \frac{dt}{t}  \Big]^{\frac pq} dr
\Big)^{1/p} < \infty\, .
\eeq
Moreover, $\| \, g \, \|^\triangle$ is equivalent to $\|\, f \, |F^s_{p,q}(\R^d)\|$.
\end{Lem}

\subsection*{The three-dimensional case}

To make use of (\ref{eq-161})  we need to know 
$\sigma_{d-1} (Q_{\lambda,t} (r))$.
For $d=3$ this quantity is explicitely known. It holds
\be
\label{pde}
\sigma_{2} (Q_{\lambda,t} (r)) = \left\{ 
\begin{array}{lll} 
\frac{\pi \, \lambda}{r} \, \Big(t^2- (\lambda - r)^2 \Big) & \qquad & \mbox{if}\quad \max (0,r-t) < \lambda < r+t
\\
&& \qquad \qquad \mbox{and}\quad t \le r+ \lambda;\\
4\pi \, \lambda^2 & &  \mbox{if}\quad t > r+\lambda;
\\
0 &&  \mbox{otherwise}\, .
\end{array}
\right.
\ee
We need a few more notation.
Depending on the behaviour of $\sigma_{2} (Q_{\lambda,t} (r)) $ we introduce the following splitting of the area of integration.
\beq\label{split1}
I_0 (r,t) & := & \{\lambda \ge 0: \quad \lambda < t-r\}\, ;
\\
I_1 (r,t) & := & \{\lambda \ge 0: \quad |\lambda-r|< \frac 34 \, t\, , 
\: \, t-r \le \lambda\}\, ;
\\
\label{split2}
I_{2} (r,t) & := & \{\lambda \ge 0: \quad  
r+ \frac 34 \, t  \le  \lambda < r+ t\, , \: \, t-r \le \lambda\}\,  ;
\\
\label{split3}
I_{3} (r,t) & := & 
\{\lambda \ge 0: \quad  
|r- t| \le   \lambda \le  r - \frac 34 \, t\}\, . 
\eeq
Obviously it follows 
\[
(\max (0,r-t), r+t) \cap [t-r,\infty)= I_1 (r,t) \cup  I_{2} (r,t) \cup  I_{3} (r,t) \, .
\]
On $I_1$ we can use 
\be\label{final1}
\frac{7}{16} \, t^2 \le \Big(t^2- (\lambda - r)^2 \Big) \le t^2\, .
\ee
Because of 
\[
\frac 74 \, t \le t+ \lambda -r \le 2\, t\, , \qquad \lambda \in I_{2}\, , 
\]
 and
\[
\frac 74 \, t \le t- \lambda + r \le 2\, t \, , \qquad \lambda \in I_{3}\, ,  
\]
we obtain 
\be\label{final2}
\frac{7}{4} \, t\, (t-\lambda +r)   \le   \Big(t^2- (\lambda - r)^2 \Big) \le  \, t\, (t-\lambda +r)\,  , 
\qquad \lambda \in I_{2}\, , 
\ee
and  
\be\label{final3}
\frac{7}{4} \, t\, (t+ \lambda -r)   \le   \Big(t^2- (\lambda - r)^2 \Big) \le  \, t\, (t+\lambda -r)\,  , 
\qquad \lambda \in I_{3}\, , 
\ee
Inserting (\ref{final1}), (\ref{final2}) and  the value of 
$\sigma_2$ in case $\lambda < t-r$ into (\ref{eq-161}), we get the following 
characterization of $RF^s_{p,q}(\R^3)$ in terms of differences.

\begin{T}\label{thm:1bb}
Let $0< p < \infty$, $0 < q \le \infty$, $1 \le v < \infty$, $0 < u \le v$  and 
\[
3 \, \max \Big(0, \frac 1p -\frac 1v, ~\frac 1q -\frac 1v \Big)< s < 1\, .
\]
Then the radial function $f \in L_p (\R^3)$ belongs to $F^s_{p,q}(\R^3)$  if, and only if,
$ g:=\tr f$ satisfies
\beqq
&& \hspace*{-0.7cm}
\| \, g \, \|^\triangle := \Big(\int_{0}^\infty r^{2}\, |g(r)|^p \, dr\Big)^{1/p}
\\
& + &
\Big(\int_{0}^1 r^{2}
\Big[\int^1_r t^{-sq} \, \Big(t^{-3} \int_{I_{0} (r,t)} \lambda^ 2 \,  |g(\lambda) - g(r)|^u\, d\lambda \Big)^{q/u} \, \frac{dt}{t}  \Big]^{p/q} \, dr
\Big)^{1/p}
\\
& + &
\Big(\int_{0}^\infty r^{2-p/u}
\Big[\int^1_0 t^{-sq} \, \Big(t^{-1} \int_{I_{1} (r,t)} \lambda \,  |g(\lambda) - g(r)|^u\, d\lambda \Big)^{q/u} \, \frac{dt}{t}  \Big]^{p/q} \, dr
\Big)^{1/p}
\\
& + &
\Big(\int_{0}^\infty r^{2-p/u}
\Big[\int^1_0 t^{-sq}  \Big(t^{-2} \, 
\int_{I_{2}(r,t)}  \lambda \, (t-\lambda +r)  \, |g(\lambda) - g(r)|^u d\lambda \Big)^{q/u}  \frac{dt}{t}  \Big]^{p/q}  dr
\Big)^{1/p} 
\\
& + &
\Big(\int_{0}^\infty r^{2-p/u}
\Big[\int^1_0 t^{-sq}  \Big(t^{-2} 
\int_{I_{3}(r,t)}  \lambda\,  (t+\lambda -r)  \, |g(\lambda) - g(r)|^u d\lambda \Big)^{q/u}  \frac{dt}{t}  \Big]^{p/q}  dr
\Big)^{1/p} < \infty  .
\eeqq
Moreover, $\| \, g \, \|^\triangle$ is equivalent to $\|\, f \, |F^s_{p,q}(\R^3)\|$.
\end{T}

\begin{Rem}
 \rm
(i) The characterization of $RF^s_{p,q}(\R^3)$ as given above seems to be of limited use because of its complexity.
From our point of view the value of Thm. \ref{thm:1bb} consists in exactly this negative  observation.
However, (\ref{eq-161}) can be used to derive an embedding into weighted spaces, see Subsection \ref{erg3} below.
\\
(ii) The formula (\ref{pde}) for $ \sigma_{2} (Q_{\lambda,t} (r))$
plays a role when calculating the solution of  the following wave equation in three dimensions
\beqq
u_{tt} (x,t) & = &  c^2 \, \Delta u (x,t)\, ,  \qquad x \in \R^3\, , \quad t >0\, ,
\\
u(x,0) & = &  0\, , \qquad  x \in \R^3\, ,
\\
u_t (x,0) & = & \left\{\begin{array}{lll}
                    1 &\qquad & \mbox{if} \quad |x| < \varrho\, ,
\\
0 && \mbox{otherwise}\, .
\end{array}\right.
\eeqq
We refer to textbooks of Strauss \cite[9.2,~Exercize~6]{Stb}
and of Dr{\'a}bek, Holubov{\'a} \cite[13.7,~Exercize~6]{DH} in this connection.
\end{Rem}

\subsection*{The two-dimensional case}

Elementary calculations, based on the law of Cosines, lead to the following formula in case $d=2$:
\be
\label{pde2}
\sigma_{1} (Q_{\lambda,t} (r)) = \left\{ 
\begin{array}{lll} 
2 \, \lambda \arccos\Big(\frac{r^2 +  \lambda^ 2 - t^ 2}{2\, r\, \lambda}\Big)  & \qquad & \mbox{if}\quad \max (0,r-t) < \lambda < r+t
\\
&& \qquad \qquad \mbox{and}\quad t \le r+ \lambda;\\
2\, \pi \, \lambda & &  \mbox{if}\quad t > r+\lambda;
\\
0 &&  \mbox{otherwise}\, .
\end{array}
\right.
\ee
According to the behaviour of $2 \, \lambda \arccos\Big(\frac{r^2 +  \lambda^ 2 - t^ 2}{2\, r\, \lambda}\Big)$
one has to split the area of integration in order to simplify (\ref{eq-161}).
We omit details.

%&&&&&&&&&&&&&&&&&&&&&&&&&&&&&&&&&&&&&&&&&&&&&&&&&&&&&&&&&&&&&&&&&&&&&&&&&&&&&&&
%&&&&&&&&&&&&&&&&&&&&&&&&&&&&&&&&&&&&&&&&&&&&&&&&&&&&&&&&&&&&&&&&&&&&&&&&&&&&&&&

\section{The characterization of radial Besov spaces  by differences}

%&&&&&&&&&&&&&&&&&&&&&&&&&&&&&&&&&&&&&&&&&&&&&&&&&&&&&&&&&&&&&&&&&&&&&&&&&&&&&&&
%&&&&&&&&&&&&&&&&&&&&&&&&&&&&&&&&&&&&&&&&&&&&&&&&&&&&&&&&&&&&&&&&&&&&&&&&&&&&&&&

Our strategy is the same as in case of Lizorkin-Triebel spaces.

%&&&&&&&&&&&&&&&&&&&&&&&&&&&&&&&&&&&&&&&&&&&&&&&&&&&&&&&&&&&&&&&&&&&&&&&&&&&&&&&
%&&&&&&&&&&&&&&&&&&&&&&&&&&&&&&&&&&&&&&&&&&&&&&&&&&&&&&&&&&&&&&&&&&&&&&&&&&&&&&&

\subsection{Differences and radial Besov spaces}

%&&&&&&&&&&&&&&&&&&&&&&&&&&&&&&&&&&&&&&&&&&&&&&&&&&&&&&&&&&&&&&&&&&&&&&&&&&&&&&&
%&&&&&&&&&&&&&&&&&&&&&&&&&&&&&&&&&&&&&&&&&&&&&&&&&&&&&&&&&&&&&&&&&&&&&&&&&&&&&&&

The point of departure is the following counterpart of Lemma \ref{prep}.

\begin{Lem}\label{prepc}
Suppose $0 <p ,q\le \infty$,  $1 \le v \le \infty$,
$0 < u \le v$,
\be\label{eq-20b}
d\, \max \Big(0\, , \frac 1p - \frac 1v  \Big)< s < 1
\ee
and  $T >0$.
Let $\Omega_t (x)$, $t>0$, $x\in \Rd$, be a family of open sets in $\Rd$ s.t.
(\ref{eq-05}) is satisfied for some $0 < A < B < \infty$.
Then $B^s_{p,q}(\Rd)$ is the collection of all
$f \in L_{\max(p,v)} (\Rd)$ s.t.
\be\label{eq-14b}
\| \, f |B^s_{p,q} (\Rd)\|^* :=
\| \, f |L_p (\Rd)\| +  \Big(\int^T_0 t^{-sq} \, \|\, (M^\Omega_{t,u} f(\, \cdot \,))|L_p (\Rd)\|^ q \, \frac{dt}{t}  \Big)^{1/q} \, .
\ee
Moreover, $\| \, \cdot \,  |B^s_{p,q} (\Rd)\|^*$ is  equivalent to
$\| \, \cdot \,  |B^s_{p,q} (\Rd)\|$
on $L_{\max(p,v)} (\Rd)$.
\end{Lem}

\noindent
\bpr
As in case of Lemma \ref{prep} the assertion represents an easy modification of Thm.~3.5.3(ii) in \cite[3.5.3]{Tr92}.
\epr

Again a Sobolev-type  embedding
\[
B^s_{p,q}(\R^d) \hookrightarrow L_u (\Rd)\, , \qquad p \le u < \left\{ \begin{array}{lll}
\frac{d}{\frac dp -s} & \qquad & \mbox{if} \quad s < d/p;
\\
\infty && \mbox{if} \quad s> d/p\, ,
\end{array}\right.
\]
see, e.g., \cite[2.7.1]{Tr83}, guarantees that $\tr$ is well-defined on $RB^s_{p,q} (\Rd)$ if $s>\sigma_{p,p} (d)$.
Now we are making use of the same arguments as in case of Lizorkin-Triebel spaces.
The radiality of $M^\Omega_{t,u} f$ and the obvious counterpart of Step 3 of the proof of  Thm. \ref{thm:1}
yield  the following characterization of radial Besov spaces.

\begin{T}\label{thm:1b}
Let $d  \ge 2$,  $0< p,q \le \infty$ and $ d/p < s < 1 $.
Then the radial function  $f \in L_p (\Rd)$ belongs to $B^s_{p,q}(\R^d)$ if, and only if,
$g :=\tr f$ satisfies
\beqq
&& \hspace*{-0.7cm}
\| \, g \, \|^\# := \|\, g\, |L_p(\R,|t|^{d-1})\|
\\
& + &
\Big(\int_0^1  t^{-sq} \Big[\int_{-\infty}^\infty |r|^{d-1} \, \Big(t^{-1} \sup_{|w|\le t} |g(r+ w) - g(r)| \Big)^{p} \, dr \Big]^{q/p}
\frac{dt}{t}  \Big]^{1/q}  < \infty \, .
\eeqq
Moreover, $\| \, g \, \|^\#$ is equivalent to $\|\, f \, |B^s_{p,q}(\R^d)\|$.
\end{T}

The same arguments, leading to Lemma \ref{diff1}, may be used to derive its counterpart for Besov spaces.

\begin{Lem}\label{diff1b}
Let $d \ge 2$, $0 < p,q \le \infty$, $1 \le v < \infty$, $0 < u \le v$  and
\[
d\, \max \Big (0\, , \frac 1p - \frac 1v \Big) < s < 1\, .
\]
Let $T >0$.
Then the radial function $f \in L_p (\Rd)$ belongs to $B^s_{p,q}(\R^d)$ if, and only if,
$g :=\tr f$ satisfies
\beq\label{eq-161b}
&& \hspace*{-0.7cm}
\| \, g \, \|^\triangle := \|\, g\, |L_p(\R,|t|^{d-1})\|
\\
\nonumber
& & \hspace{-0.7cm}
\Big(\int^T_0 t^{-sq}
\Big[\int_{0}^\infty r^{d-1}   \Big(t^{-d} \int_{\max(0,r-t)}^{r+t}     |g(\lambda) - g(r)|^u
\sigma_{d-1} (Q_{\lambda,t} (r) )  \, d\lambda \Big)^{\frac pu} dr \Big]^{q/p} \frac{dt}{t}  \Big)^{1/q}
< \infty\, .
\eeq
Moreover, $\| \, g \, \|^\triangle$ is equivalent to $\|\, f \, |B^s_{p,q}(\R^d)\|$.
\end{Lem}

By restricting to  $d=3$ and by employing the splitting from (\ref{split1})-(\ref{split3}) this results in the following.

\begin{T}\label{thm:1bbb}
Let $0 <p, q \le \infty$, $1 \le v < \infty$, $0 < u \le v$  and 
\[
3 \, \max \Big(0, \frac 1p -\frac 1v\Big)< s < 1\, .
\]
Then the radial function $f \in L_p (\R^3)$ belongs to $B^s_{p,q}(\R^3)$  if, and only if,
$ g:=\tr f$ satisfies
\beqq
&& \hspace*{-1.0cm}
\| \, g \, \|^\triangle := \Big(\int_{0}^\infty r^{2}\, |g(r)|^p \, dr\Big)^{1/p}
\\
& + &
\Big(\int^1_0 t^{-sq}
\Big[\int_{0}^t r^{2}
\Big(t^{-3} \int_{I_{0} (r,t)} \lambda^2   |g(\lambda) - g(r)|^u d\lambda \Big)^{p/u}  dr
\Big]^{q/p} \frac{dt}{t}  \Big)^{1/q} 
\\
& + &
\Big(\int^1_0 t^{-sq}
\Big[\int_{0}^\infty r^{2-p/u}
\Big(t^{-1} \int_{I_{1} (r,t)} \lambda   |g(\lambda) - g(r)|^u d\lambda \Big)^{p/u}  dr \Big]^{q/p} \frac{dt}{t}  \Big)^{1/q} 
\\
& + &
\Big(\int^1_0 t^{-sq} 
\Big[ \int_{0}^\infty r^{2-p/u} \Big(t^{-2}  
\int_{I_{2}(r,t)}  \lambda \, (t-\lambda +r)  |g(\lambda) - g(r)|^u d\lambda \Big)^{p/u} dr \Big]^{q/p}  \frac{dt}{t}  \Big)^{1/q} 
\\
& + &
\Big(\int^1_0 t^{-sq} \Big[\int_{0}^\infty r^{2-p/u}  \Big(t^{-2} 
\int_{I_{3}(r,t)}  \lambda (t+\lambda -r)  |g(\lambda) - g(r)|^u d\lambda \Big)^{p/u} dr \Big]^{q/p} \frac{dt}{t}  \Big)^{1/q} < \infty  .
\eeqq
Moreover, $\| \, g \, \|^\triangle$ is equivalent to $\|\, f \, |B^s_{p,q}(\R^3)\|$.
\end{T}

%&&&&&&&&&&&&&&&&&&&&&&&&&&&&&&&&&&&&&&&&&&&&&&&&&&&&&&&&&&&&&&&&&&&&&&&&&&&&&&&
%&&&&&&&&&&&&&&&&&&&&&&&&&&&&&&&&&&&&&&&&&&&&&&&&&&&&&&&&&&&&&&&&&&&&&&&&&&&&&&&

\subsection{Some concluding remarks}
\label{erg}

%&&&&&&&&&&&&&&&&&&&&&&&&&&&&&&&&&&&&&&&&&&&&&&&&&&&&&&&&&&&&&&&&&&&&&&&&&&&&&&&
%&&&&&&&&&&&&&&&&&&&&&&&&&&&&&&&&&&&&&&&&&&&&&&&&&&&&&&&&&&&&&&&&&&&&&&&&&&&&&&&

Here we shall indicate the main problems connected
with a desirable extension of Thm.  \ref{thm:1} and Thm. \ref{thm:1b}, respectively.
\\
Smoothness $s>1$ requires the use of higher order differences.
Let $N \in \N$ and $x \in \Rd$. Then the $N$-th order difference of $f$ in $x$
is defined to be
\[
\Delta_h^N f (x) := \sum_{j=0}^N {N \choose j}\, (-1)^{N-j} \, f(x+jh)
\]
By using these higher order differences the following extension of Lemma
\ref{prep} can be proved.

\begin{Lem}\label{prepb}
Suppose $0 <p <\infty$, $0 < q \le \infty$, $1 \le v \le \infty$,
$0 < u \le v$, $N \in \N$,
\[
d\, \max \Big(0\, , \frac 1p - \frac 1v \, ,  \frac 1q - \frac 1v  \Big)< s <  N\, ,
\]
and  $T >0$.
Let $\Omega_t (x)$, $t>0$, $x\in \Rd$, be a family of open sets in $\Rd$ s.t.
(\ref{eq-05}) is satisfied for some $0 < A < B < \infty$.
Then $F^s_{p,q}(\Rd)$ is the collection of all
$f \in L_{\max(p,v)} (\Rd)$ s.t.
\[
\| \, f |F^s_{p,q} (\Rd)\|^* :=
\| \, f |L_p (\Rd)\| + \Big\| \Big(\int^T_0 t^{-sq} \, (M^{N,\Omega}_{t,u} f(\, \cdot \,))^ q \,
\frac{dt}{t}  \Big)^{1/q} \, \Big|L_p (\Rd)\Big\|\, ,
\]
where
\[
M^{N,\Omega}_{t,u} f (x):= \Big(t^{-d}\, \int_{\Omega_t (x)} | \Delta_h^N f(x)|^u \, dh \Big)^{1/u} \,  .
\]
Moreover, $\| \, \cdot \,  |F^s_{p,q} (\Rd)\|^*$ is  equivalent to
$\| \, \cdot \,  |F^s_{p,q} (\Rd)\|$
on $L_{\max(p,v)} (\Rd)$.
\end{Lem}

Of course, there is a counterpart for Besov spaces as well.

\begin{Lem}\label{prepcb}
Suppose $N \in \N$ $0 <p ,q\le \infty$,  $1 \le v \le \infty$,
$0 < u \le v$,
\[
d\, \max \Big(0\, , \frac 1p - \frac 1v  \Big)< s < N\, , 
\]
and  $T >0$.
Let $\Omega_t (x)$, $t>0$, $x\in \Rd$, be a family of open sets in $\Rd$ s.t.
(\ref{eq-05}) is satisfied for some $0 < A < B < \infty$.
Then $B^s_{p,q}(\Rd)$ is the collection of all
$f \in L_{\max(p,v)} (\Rd)$ s.t.
\[
\| \, f |B^s_{p,q} (\Rd)\|^* :=
\| \, f |L_p (\Rd)\| +  \Big(\int^T_0 t^{-sq} \, \|\, (M^{N,\Omega}_{t,u} f(\, \cdot \,))|L_p (\Rd)\|^ q \, \frac{dt}{t}  \Big)^{1/q} \, .
\]
Moreover, $\| \, \cdot \,  |B^s_{p,q} (\Rd)\|^*$ is  equivalent to
$\| \, \cdot \,  |B^s_{p,q} (\Rd)\|$
on $L_{\max(p,v)} (\Rd)$.
\end{Lem}

Now we are in position to describe the difficulties with the extension of Thm.  \ref{thm:1} and Thm. \ref{thm:1b}, respectively.
Let $f$ be a radial function. Then, as above in Steps 1,2 of the proof of Thm. \ref{thm:1}, one can show that
$M^{N,\Omega}_{t,u} f (x)$ is a radial function.
However, as it is easy to see, for $N \ge 2$  we do not  have an identity of the form
\[
\Delta_h^N f (x) = \Delta_w^N g (|x|)
\]
in general
(here $f$, $h$ and $x$ are given and we may choose $w$).
For this reason we do not have a counterpart of Step 3
(proof of Thm. \ref{thm:1}) in case of higher order differences.

%&&&&&&&&&&&&&&&&&&&&&&&&&&&&&&&&&&&&&&&&&&&&&&&&&&&&&&&&&&&&&&&&&&&&&&&&&&&&&&&
%&&&&&&&&&&&&&&&&&&&&&&&&&&&&&&&&&&&&&&&&&&&&&&&&&&&&&&&&&&&&&&&&&&&&&&&&&&&&&&&

\section{A comparison with weighted Besov-Lizorkin-Triebel spaces}
\label{sec3}

%&&&&&&&&&&&&&&&&&&&&&&&&&&&&&&&&&&&&&&&&&&&&&&&&&&&&&&&&&&&&&&&&&&&&&&&&&&&&&&&
%&&&&&&&&&&&&&&&&&&&&&&&&&&&&&&&&&&&&&&&&&&&&&&&&&&&&&&&&&&&&&&&&&&&&&&&&&&&&&&&

This section  is devoted to the comparison of $RF^s_{p,q} (\Rd)$ and $RB^s_{p,q} (\Rd)$ 
with certain weighted spaces.

%&&&&&&&&&&&&&&&&&&&&&&&&&&&&&&&&&&&&&&&&&&&&&&&&&&&&&&&&&&&&&&&&&&&&&&&&&&&&&&&
%&&&&&&&&&&&&&&&&&&&&&&&&&&&&&&&&&&&&&&&&&&&&&&&&&&&&&&&&&&&&&&&&&&&&&&&&&&&&&&&

\subsection{Besov-Lizorkin-Triebel spaces and Muckenhoupt weights}
\label{erg4}

%&&&&&&&&&&&&&&&&&&&&&&&&&&&&&&&&&&&&&&&&&&&&&&&&&&&&&&&&&&&&&&&&&&&&&&&&&&&&&&&
%&&&&&&&&&&&&&&&&&&&&&&&&&&&&&&&&&&&&&&&&&&&&&&&&&&&&&&&&&&&&&&&&&&&&&&&&&&&&&&&

For the Muckenhoupt class $\ca_\infty$ of weights  there is a well-established theory for 
Besov and Lizorkin-Triebel spaces (for definitions we refer to the Appendix),
we refer to 
Bui et all \cite{bui-1,bui-2,BPT-1,BPT-2}, Kokilashvili \cite{Ko,Ko2}, Rychkov \cite{Ry},  
Bownik et all  \cite{Bow1,Bow2}, Haroske, Piotrowska \cite{HP}, Haroske, Skrzypczak \cite{HS},
and Izuki, Sawano \cite{IS}. Concerning Muckenhoupt weights we refer to the monograph \cite{stein} of Stein. 
Recall, the function $w(t):= t^{d-1}$, $t\in \R$, is a Muckenhoupt weight for any $d \in \N$.
As a consequence of Thm. 8 and Thm. 9  in \cite{SSV1} one knows the following.

\begin{Prop}
Let $d\ge 2$, $0 < p < \infty$ and  $0 < q \le  \infty$.
\\
(i)  Let $s> \sigma_{1,q}(d)$
and either
\[
 s > d\, \Big(\frac 1p - \frac 1d\Big)
\qquad \mbox{or}\qquad
s = d\, \Big(\frac 1p - \frac 1d\Big)\qquad \mbox{and}\qquad 0 < p \le 1\, .
\]
Then we have coincidence
\[
R{F}^s_{p,q} (\R, |t|^{d-1}) = \tr (RF^s_{p,q} (\Rd)) \qquad \mbox{(in the sense of equivalent quasi-norms).}
\]
(ii) Let either
\[
 s > d\, \Big(\frac 1p - \frac 1d\Big)
\qquad \mbox{or}\qquad
s = d\, \Big(\frac 1p - \frac 1d\Big)\qquad \mbox{and}\qquad 0 < q \le 1\, .
\]
Then we have coincidence
\[
R{B}^s_{p,q} (\R, |t|^{d-1}) = \tr (RB^s_{p,q} (\Rd)) \qquad \mbox{(in the sense of equivalent quasi-norms).}
\]
\end{Prop}

Hence, we may turn our characterizations of 
$\tr (RF^s_{p,q} (\Rd))$ into characterizations of  
$R{F}^s_{p,q} (\R, |t|^{d-1})$. As an example we reformulate Thm. \ref{thm:1bb}.

\begin{Cor}\label{cor:1bb}
Let $3/2 < p < \infty$, $0 < q \le \infty$, $1 \le v < \infty$, $0 < u \le v$  and 
\[
3 \, \max \Big(0, \frac 1p -\frac 13, \frac 1p -\frac 1v, ~\frac 1q -\frac 1v \Big)< s < 1\, .
\]
Then the even  function $g \in L_p (\R, |t|^2)$ belongs to $F^s_{p,q}(\R, |t|^2)$  if, and only if,
$ g$ satisfies
\beqq
&& \hspace*{-1.0cm}
\| \, g \, \|^\triangle := \Big(\int_{0}^\infty r^{2}\, |g(r)|^p \, dr\Big)^{1/p}
\\
& + &
\Big(\int_{0}^1 r^{2}
\Big[\int^1_r t^{-sq} \, \Big(t^{-3} \int_{I_{0} (r,t)} \lambda^ 2 \,  |g(\lambda) - g(r)|^u\, d\lambda \Big)^{q/u} \, \frac{dt}{t}  \Big]^{p/q} \, dr
\Big)^{1/p}
\\
& + &
\Big(\int_{0}^\infty r^{2-p/u}
\Big[\int^1_0 t^{-sq}  \Big(t^{-1} \int_{I_{1} (r,t)} \lambda   |g(\lambda) - g(r)|^u d\lambda \Big)^{q/u}  \frac{dt}{t}  \Big]^{p/q}  dr
\Big)^{1/p}
\\
& + &
\Big(\int_{0}^\infty r^{2-p/u} \Big[\int^1_0 t^{-sq}  \Big(t^{-2} 
\int_{I_{2}(r,t)}  \lambda  (t-\lambda +r)   |g(\lambda) - g(r)|^u d\lambda \Big)^{q/u}  \frac{dt}{t}  \Big]^{p/q}  dr
\Big)^{1/p} 
\\
& + &
\Big(\int_{0}^\infty r^{2-p/u} \Big[\int^1_0 t^{-sq}  \Big(t^{-2} 
\int_{I_{3}(r,t)}  \lambda  (t+\lambda -r)   |g(\lambda) - g(r)|^u d\lambda \Big)^{q/u}  \frac{dt}{t}  \Big]^{p/q}  dr
\Big)^{1/p} < \infty  .
\eeqq
Moreover, $\| \, g \, \|^\triangle$ is equivalent to $\|\, g \, |F^s_{p,q}(\R,|t|^2)\|$.
\end{Cor}

\begin{Rem}
 \rm
Cor. \ref{cor:1bb} may be understood as a first hint how complicated characterizations of 
weighted spaces ${F}^s_{p,q} (\R,w)$ with $w$ being a Muckenhoupt weight  may look like.
This is in certain contrast to the 
characterization of   ${F}^s_{p,q} (\R,w)$ with smooth weights  $w$, see 
\cite[5.1]{ST} and the next subsection.
\end{Rem}

%&&&&&&&&&&&&&&&&&&&&&&&&&&&&&&&&&&&&&&&&&&&&&&&&&&&&&&&&&&&&&&&&&&&&&&&&&&&&&&&
%&&&&&&&&&&&&&&&&&&&&&&&&&&&&&&&&&&&&&&&&&&&&&&&&&&&&&&&&&&&&&&&&&&&&&&&&&&&&&&&

\subsection{Weighted Besov-Lizorkin-Triebel spaces and differences}
\label{erg3}

%&&&&&&&&&&&&&&&&&&&&&&&&&&&&&&&&&&&&&&&&&&&&&&&&&&&&&&&&&&&&&&&&&&&&&&&&&&&&&&&
%&&&&&&&&&&&&&&&&&&&&&&&&&&&&&&&&&&&&&&&&&&&&&&&&&&&&&&&&&&&&&&&&&&&&&&&&&&&&&&&

In the monograph \cite{ST} Schmei{\ss}er and Triebel developed the theory of 
weighted Besov-Lizorkin-Triebel spaces for a certain class of smooth weights.
This class of weights cover the following expamples:
\[
{\rho}_\alpha (x) := (1+ |x|^2)^{\alpha/2}, \qquad x \in \Rd\, , \qquad \alpha >0\,,
\]
see Remark 2 in \cite[1.4.1]{ST}. We concentrate on $\alpha = d-1$, $d\in \N$.
In this situation Thm.~5.1.4 in \cite{ST} reads as follows.

\begin{Prop}\label{diffi}
Let $d\in \N$, $d\ge 2$.\\
 (i) Let $0 < p,q < \infty$, let $N \in \N$ and $\sigma_{p,q}(1) < s < N$.
Then $\| \, f \, |{F}^s_{p,q} (\R,\rho_{d-1})\|$ and
\beqq
%&& \hspace*{-1.0cm}
\| \, f \, |{F}^s_{p,q} (\R, \rho_{d-1})\|^* & := &  \| \, f  \,  |L_p (\R,\rho_{d-1})\|
\\
& + &
\Big\| \,  \Big( \int_0^1  t^{-sq} \, \Big[ \int_{|h| \le 1} |\Delta_{ht}^N f (\, \cdot \, )|\, dh \Big]^{q}
\, \frac{dt}{t}  \Big)^{1/q} \, \Big|L_p (\R, \rho_{d-1})\Big\|
\eeqq
are equivalent.
\\
(ii) Let $0 < p <  \infty$, $0 < q \le \infty$,  $N \in \N$ and $\sigma_{p,p}(1) < s < N$.
Then $\| \, f \, | {B}^s_{p,q} (\R,\rho_{d-1})\|$ and
\beqq
%&& \hspace*{-1.0cm}
\| \, f \, | {B}^s_{p,q} (\R, \rho_{d-1})\|^* := \| \, f  \,  |L_p (\R, \rho_{d-1})\|
+
\Big( \int_{|h|\le 1}  |h|^{-sq} \, \| \, \Delta_{h}^N f\,  |L_p (\R, \rho_{d-1})\|^{q}
\, \frac{dh}{|h|^d} \Big)^{1/q}
\eeqq
are equivalent.
\end{Prop}

Our next step consists in taking the characterization in Prop. \ref{diffi} and 
comparing this  with formula (\ref{eq-161}) in case $u=1$.
Because of
\[
\sigma_{d-1} (Q_{\lambda,t} (r)) \le \sigma_{d-1} (\{x: \: |x|=t\}) \lsim t^{d-1}
\]
and 
\beqq
&& \hspace*{-0.8cm}\int_{0}^\infty x^{d-1}
\Big[\int^1_0 t^{-sq} \, \Big(t^{-1} \int_{|h|\le t} |g(x+ h) - g(x)|\, dh \Big)^{q} \, \frac{dt}{t}  \Big]^{p/q} \, dx
\\
& \lsim &
\int_{-\infty}^\infty \, (1+ |x|^2)^{(d-1)/2} \,
\Big[ \int_0^1  t^{-sq} \, \Big( t^{-1} \, \int_{|h| \le t} |g (x+h) - g(x)|\, dh \Big)^{q}
\, \frac{dt}{t}  \Big]^{p/q} \, dx
\eeqq
we conclude the continuous embedding
\[
R{F}^s_{p,q} (\R, \rho_{d-1}) \hookrightarrow
\tr (RF^s_{p,q} (\Rd))\, ,
\]
if $0 < p,q < \infty$ and $\sigma_{p,q}(d) < s < 1$. The same type of arguments leads to
\[
R{B}^s_{p,q} (\R, \rho_{d-1}) \hookrightarrow
\tr (RB^s_{p,q} (\Rd))\, ,
\]
if $0 < p,q \le \infty$ and $\sigma_{p,p}(d) < s < 1$.
\\
Both assertions supplement Prop. \ref{diffi} in view of the trivial 
embeddings
\[
{F}^s_{p,q} (\R, \rho_{d-1}) \hookrightarrow
F^s_{p,q} (\R, |t|^{d-1})
\qquad \mbox{and}\qquad {B}^s_{p,q} (\R, \rho_{d-1}) \hookrightarrow
B^s_{p,q} (\R, |t|^{d-1})
\]
(valid for all admissible parameters as a direct consequence of the definition).

%&&&&&&&&&&&&&&&&&&&&&&&&&&&&&&&&&&&&&&&&&&&&&&&&&&&&&&&&&&&&&&&&&&&&6
%&&&&&&&&&&&&&&&&&&&&&&&&&&&&&&&&&&&&&&&&&&&&&&&&&&&&&&&&&&&&&&&&&&&&&

\section{Appendix -- Muckenhoupt weights and function spaces}

%&&&&&&&&&&&&&&&&&&&&&&&&&&&&&&&&&&&&&&&&&&&&&&&&&&&&&&&&&&&&&&&&&&&&6
%&&&&&&&&&&&&&&&&&&&&&&&&&&&&&&&&&&&&&&&&&&&&&&&&&&&&&&&&&&&&&&&&&&&&&

For convenience of the reader we collect some definitions and properties
around Muckenhoupt weights and associated weighted function spaces.

%&&&&&&&&&&&&&&&&&&&&&&&&&&&&&&&&&&&&&&&&&&&&&&&&&&&&&&&&&&&&&&&&&&&&6
%&&&&&&&&&&&&&&&&&&&&&&&&&&&&&&&&&&&&&&&&&&&&&&&&&&&&&&&&&&&&&&&&&&&&&

\subsection{Muckenhoupt weights}

%&&&&&&&&&&&&&&&&&&&&&&&&&&&&&&&&&&&&&&&&&&&&&&&&&&&&&&&&&&&&&&&&&&&&6
%&&&&&&&&&&&&&&&&&&&&&&&&&&&&&&&&&&&&&&&&&&&&&&&&&&&&&&&&&&&&&&&&&&&&&

A weight function (or simply a weight) is a nonnegative, locally integrable function on $\Rd$.
We collect a few facts including the definition of Muckenhoupt weights.
As usual, $p'$ is related to $p$ via the formula  $1/p + 1/p'=1$.

\begin{Def}
Let $1 < p < \infty$.
Let $w$ be a nonnegative, locally integrable function on $\Rd$.
Then $w$ belongs to the  Muckenhoupt class $\ca_p$, if there exists a constant $A>0$ s.t. for all balls $B$ the following inequality holds:
\[
\Big(\frac{1}{|B|} \int_B w(x)\, dx\Big)^{1/p} \, \cdot \,
\Big(\frac{1}{|B|} \int_B w(x)^{-p'/p} \, dx\Big)^{1/p'} \le A\, .
\]
\end{Def}

The Muckenhoupt class $\ca_\infty$ is defined as
\[
\ca_\infty := \bigcup_{p>1} \ca_p\, .
\]

\begin{Rem}
 \rm
(i) Good sources for Muckenhoupt weights are Stein's monograph \cite{stein} and the graduate text \cite{Du} of Duoandikoetxea.
\\
(ii) It is well-known that the functions 
\[
w(t):= |t|^{d-1}, \quad t\in \R, \qquad  \mbox{and} \qquad
\rho_{d-1} (t):= (1+|t|^2)^{(d-1)/2}, \quad t\in\R, 
\]
belong to $\ca_\infty$.
\end{Rem}

%&&&&&&&&&&&&&&&&&&&&&&&&&&&&&&&&&&&&&&&&&&&&&&&&&&&&&&&&&&&&&&&&&&&&6
%&&&&&&&&&&&&&&&&&&&&&&&&&&&&&&&&&&&&&&&&&&&&&&&&&&&&&&&&&&&&&&&&&&&&&

\subsection{Weighted Besov and Lizorkin-Triebel spaces}

%&&&&&&&&&&&&&&&&&&&&&&&&&&&&&&&&&&&&&&&&&&&&&&&&&&&&&&&&&&&&&&&&&&&&6
%&&&&&&&&&&&&&&&&&&&&&&&&&&&&&&&&&&&&&&&&&&&&&&&&&&&&&&&&&&&&&&&&&&&&&

To introduce weighted  Besov and Lizorkin-Triebel spaces we make use of tools from Fourier analysis.
Let $\psi \in C_0^\infty (\Rd)$ be a  function such that
$\psi (x)=1$, $|x| \le 1$ and $\psi (x)=0$, $|x| \ge 2$. Then
by
\be\label{ddu}
\varphi_0 (x) := \psi (x), \quad \varphi_1 (x):=
\varphi_0 (x/2) - \varphi_0 (x), \quad \varphi_j (x) := \varphi_1 (2^{-j+1}x),
\ee
$j \in \N$, we define a smooth dyadic decomposition of unity.
By using such a decomposition of unity we introduce  weighted function spaces as follows.

\begin{Def}
Let $0 < q \le \infty$, $s\in \R$ and $w \in \ca_\infty$.\\
{\rm (i)}
Let $0 < p< \infty$. Then the weighted Besov space
$B^s_{p,q}(\Rd,w)$ is the collection of all $f \in S'(\Rd)$
such that
\[
\| \, f \, |B^s_{p,q}(\Rd,w)\| := \Big(
\sum_{j=0}^\infty 2^{jsq}\, \| \, \cfi [\varphi_j (\xi)\, \cf f (\xi)](\, \cdot  \, )\, |L_p(\Rd,w)\|^q\Big)^{1/q}<\infty\, .
\]
{\rm (ii)}
Let $0 < p< \infty$. Then the weighted Triebel-Lizorkin space
$F^s_{p,q}(\Rd,w)$ is the collection of all $f \in S'(\Rd)$
such that
\[
\| \, f \, |F^s_{p,q}(\Rd,w)\| := \Big\|\, \Big(
\sum_{j=0}^\infty 2^{jsq}\, | \, \cfi [\varphi_j (\xi)\, \cf f (\xi)](\, \cdot  \, )\, |^q \Big)^{1/q} \Big| L_p (\Rd,w)\Big\| <\infty\, .
\]
{\rm (iii)}
Let $p= \infty$. Then the  Besov space
$B^s_{\infty,q}(\Rd)$ is the collection of all $f \in S'(\Rd)$
such that
\[
\| \, f \, |B^s_{\infty,q}(\Rd)\| := \Big(
\sum_{j=0}^\infty 2^{jsq}\, \| \, \cfi [\varphi_j (\xi)\, \cf f (\xi)](\, \cdot  \, )\, |L_\infty(\Rd)\|^q\Big)^{1/q}<\infty\, .
\]
\end{Def}

\begin{Rem}
\label{unendlich}
\rm
(i)
For $w\equiv 1$ we are in the unweighted case.
The associated spaces are denoted by $B^s_{p,q}(\Rd)$ and $F^s_{p,q}(\Rd)$.
\\
(ii)
Observe, that we did not define weighted spaces with $p=\infty$. However, it will be convenient for us to use
the convention $B^s_{\infty,q}(\Rd,w):= B^s_{\infty,q}(\Rd)$.
\\
(iii) Weighted Besov and Lizorkin-Triebel spaces with $w\in {\ca}_\infty$ have been first studied
systematically by Bui \cite{bui-1,bui-2}, cf. also \cite{BPT-1}  and \cite{BPT-2}.
In addition we refer  to Kokilashvili \cite{Ko,Ko2},
Haroske, Piotrowska \cite{HP} and \cite{HS}.
Standard references for unweighted spaces are
the monograph's \cite{Pe,Tr83,Tr92,Tr06} as well as \cite{FJ2}.
More general classes of weights have been treated by Rychkov \cite{Ry}, 
Bownik and Ho \cite{Bow1}, Bownik \cite{Bow2},  Izuki and Sawano \cite{IS}, and Wojciechowska \cite{AW,AW2}.
\end{Rem}

%&&&&&&&&&&&&&&&&&&&&&&&&&&&&&&&&&&&&&&&&&&&&&&&&&&&&&&&&&&&
%&&&&&&&&&&&&&&&&&&&&&&&&&&&&&&&&&&&&&&&&&&&&&&&&&&&&&&&&&&&

\end{document}